\documentclass[preprint,12pt]{elsarticle}

\usepackage{amssymb}
\usepackage{amsmath}
\usepackage{subfigure}

\journal{Journal of Computational Physics}

\begin{document}

\begin{frontmatter}

\title{Lagrangian Particle Method for Compressible Fluid Dynamics}

\author[label1]{Hsin-Chiang Chen}
\author[label1,label2]{Roman Samulyak\footnote{Corresponding author, roman.samulyak@stonybrook.edu}}
\author[label1]{Wei Li}

\address[label1]{Department of Applied Mathematics and Statistics,
  Stony Brook University, Stony Brook, NY 11794}
\address[label2]{Computational Science Center, Brookhaven National
  Laboratory, Upton, NY 11973}

\begin{abstract}

A new Lagrangian particle method for solving Euler
equations for compressible inviscid  fluid or gas flows is proposed.  
Similar to smoothed particle hydrodynamics (SPH),
the method represents fluid cells with Lagrangian particles and is suitable for the
simulation of complex free surface / multiphase flows.
The main contributions of our method, which is different from SPH in all other aspects,
are (a) significant improvement of approximation of differential operators based on 
a polynomial fit via  weighted least squares approximation
and the convergence of prescribed order, (b) an upwinding second-order particle-based algorithm with limiter,
providing  accuracy and long term stability, (c)  elimination of the dependence on artificial parameters such as the
smoothening length in SPH, causing difficulties in the case
of large density changes, and (d) accurate resolution of states at free interfaces.
Numerical verification test demonstrating the convergence order are presented as well as examples of complex 
free surface flows.

\end{abstract}

\begin{keyword}
Lagrangian fluid mechanics \sep particle method \sep generalized
finite differences

\MSC 65M06 \sep 70F99 \sep 76T10

\end{keyword}
\end{frontmatter}

\section{Introduction and Motivation}
\label{sec:introduction}

High resolution Lagrangian methods are essential for achieving
predictive simulations of a wide spectrum of complex free surface /
multiphase problems. Most widely used approaches for the simulation of
multiphase problems are based on Eulerian meshes enhanced with special
algorithms for resolving interfaces such as the volume-of-fluid
\cite{Hirt81}, the level set method \cite{Osher02}, arbitrary
Lagrangian -€" Eulerian methods \cite{ALE}, or the method of front
tracking \cite{FTcode} which is a hybrid method involving a moving
Lagrangian mesh over a fixed Eulerian mesh. In addition, they often
use various adaptive features such as adaptive mesh refinement.
These and finite element methods, most common for engineering problems with irregular
geometries, require complex computationally intensive methods for the
generation of high quality meshes.

Theoretically, the traditional Lagrangian formulation of fluid
dynamics \cite{Lamb} is the basis for the most natural and accurate
method for the simulation of complex free surface and multiphase
systems, but it suffers from the mesh distortion problem in unsteady,
turbulent flows. As a result, the Lagrangian methods are widely used
only in 1D for all problems except the dynamics of solids that is
characterized by small deformations. The overwhelming majority of
solid dynamics codes use finite element-based Lagrangian methods,€" the
fact that speaks for advantages of Lagrangian approaches within their
applicability range.

A way to extend the Lagrangian method to 3D was proposed in smooth
particle hydrodynamics (SPH). SPH \cite{Monaghan92,Monaghan05}
is a Lagrangian particle method
in computational fluid dynamics in which deforming Lagrangian cells
are replaced with particles. SPH eliminates the main mesh tangling
difficulty of the original Lagrangian method while retaining many of
its advantages. Due to its
Lagrangian nature, SPH is strictly mass-conservative and capable of
robustly handling interfaces of arbitrary complexity in the simulation
of free surface and multiphase flows. The representation of matter by
particles provides adaptivity to density changes. Not only does this
improve the traditional adaptive mesh refinement of structural meshes,
that introduces sharp boundaries between mesh patches of different
resolution, but also it enables simulations of large ranges of spatial
scales (for instance expansion into vacuum and matter islands
separated by large vacuum domains).

However the major drawback of SPH is a very poor accuracy of discrete
differential operators. It is widely accepted \cite{Diltz99,GIZMO},
including original SPH developers \cite{Monaghan05}, that the
traditional SPH discretization has zero-order convergence for widely
used kernels. The SPH discretization of derivatives is convergent to 
the order consistent with the interpolating polynomial for the kernel
function only if particles are located on a rectangular mesh (which 
is not the case for unsteady flows). 
In addition, it depends on artificial parameters, in
particular on the smoothening radius, causing major difficulties in
the case of large density changes. The reason why SPH produces stable
and reasonable results for certain problems, despite using inaccurate
and non-convergent discretization of differential operators, is its
connection to the Lagrangian / Hamiltonian dynamics of particles \cite{Price_MHD}.
In particular, the traditional discrete SPH equations for the compressible Euler equations
are not accurate, but they
accurately represent equations of the Lagrangian dynamics of particles interacting via
isentropic potentials. The Lagrangian and Hamiltonian properties are also 
responsible for the long term stability of the traditional SPH.
But the Hamiltonian dynamics of particles only approximately
represent the dynamics of continuum hydrodynamic systems, and the 
isentropic interaction energy places additional restrictions.
 
A number of 'modern' or 'corrected' SPH methods
have been developed in recent years (see \cite{Diltz99} and reviews 
\cite{Monaghan05,GIZMO}). They
include the moving-least-squares SPH, 'Godunov'-SPH, P-SPH, PHANTOM
etc. But they all improve certain features of SPH at the expense of other properties such as
conservation, long-time stability, or prohibitively large number of
neighbors that causes other problems. They all still have zero-convergence order,
except for the 1st order convergent, moving-least-squates SPH \cite{Diltz99},
that suffers from long-time stability and other issues. 
A new class of convergent, mesh-free hydrodynamic simulation methods was developed in \cite{GIZMO}.

We have proposed a new Lagrangian particle method for solving
compressible Euler equations that eliminates major deficiencies of
SPH: the dependence on a parameter called smoothening length, the
presence of large linear errors in SPH differential operators, and 
ensures long term stability via upwinding discretization methods.
Significantly different from SPH in most of approximations, our method is also
easily generalizable to coupled system of hyperbolic and elliptic
or parabolic PDE's for other physics processes.

In the Lagrangian particle method, approximations of spatial
derivatives are obtained by employing a local polynomial fit known also as
the generalized finite difference (GFD) method \cite{BenitoUrena01}.  The main idea is
to find closest neighbors of each particle and approximate the spatial
derivative of a certain physical quantity the particle location as a
linear combination of this quantity at neighboring particles. The
optimal coefficients in this linear combination are calculated by
solving a least squares problem. Second order accurate spatial
discretization is used in the current algorithm, but the GFD method
makes it possible to use higher order discretizations with increased
particle neighborhoods.  Our algorithms uses much smaller number of neighbor particles
compared to the Godunov-SPH or other recent SPH modifications
that may require hundreds of particles \cite{GIZMO}. An application example 
of the GFD method to the advection-diffusion equation is given in \cite{GFD_ad}.

The conservative Lagrangian  formulation of
Euler equations is transformed into a quasi-linear form, and an
upwinding scheme is employed for the numerical integration. 
Multiple spatial dimensions are resolved using a Strang splitting method for 
Euler equations.  Research on algorithms for elliptic problems
involving geometrically complex boundaries and interfaces is in
progress.  Together with hyperbolic solvers, they form the basis for
simulations of complex multiphysics and multiphase systems.  The
Lagrangian particle method has been implemented in all dimensions, but 
algorithms for a parallel code for solving three-dimensional Euler equations 
will be reported in a forthcoming paper. The main goal of the present work is 
to introduce main ideas of the Lagrangian particles dynamics and to present verification tests.  

The paper is organized as follows. Section 2 contains the main
governing equations.  Numerical discretization and main algorithms of
the Lagrangian particle dynamics are presented in Section 3. Section 4
presents verification tests and accuracy studies.  We conclude this work by a 
summary of our results and perspectives for the future work.

\section{Governing Equations}

Consider the one-dimensional Lagrangian formulation of the Euler
equations, written in the conservative form
\cite{RichtmyerMorton,Lamb}
\begin{equation} \label{euler1}
	U_{t}^{'} + \left[ F(U^{'}) \right]_{x} = 0,
\end{equation}
\begin{equation} \label{euler2}
	U^{'} = \left( \begin{array} {c} V \\ u \\ E \end{array} \right), 
	\, F(U^{'}) = V_{0} \left( \begin{array} {c} -u \\ P \\ Pu \end{array} \right),
\end{equation}
where $V$ is the specific volume, $u$ is the velocity, $E$ is the specific
total energy, and $P$ is the pressure.

Let's assume that the equation of state (EOS) is in the form 
$e=f\left( P,V \right) $, 
where $e$ is the specific internal energy, $e = E - u^2/2$.
Equations (\ref{euler1}) and (\ref{euler2}) can be written using 
$U=[V$ $u$ $P]^{T}$ as the state vector as follows
\begin{equation} \label{euler3}
	U_{t} + A(U)U_{x} = 0,
\end{equation}
\begin{equation} \label{euler4}
	U = \left( \begin{array} {c} V \\ u \\ P \end{array} \right), \quad
	A(U) = V_{0}\left( \begin{array} {ccc} 
	0 & -1 & 0 \\ 
	0 & 0 & 1  \\ 
	0 & K & 0 \end{array} \right),
\end{equation}
where 
\begin{equation} \label{Kdef}
	K = \left(P+\frac{\partial e}{\partial V}\right)\left/\frac{\partial e}{\partial P}\right. .
\end{equation}
For example, using the polytropic gas EOS
\begin{equation} \label{polyEOS}
	e = \frac{PV}{\gamma -1},
\end{equation} 
where $\gamma$ is the ratio of specific heats, we obtain 
\begin{equation}
	A(U) = V_{0}\left( \begin{array} {ccc} 
	0 & -1 & 0 \\ 
	0 & 0 & 1  \\ 
	0 & (\frac{c}{V})^{2} & 0 \end{array} \right),
\end{equation}
where $c=\sqrt{\gamma PV}$ is the speed of sound. 
Note that the transformation is exact 
(i.e. not a result of linerization).
If the matrix $A$ is diagonalized as $A = R\Lambda R^{-1}$, 
equations (\ref{euler3}) and (\ref{euler4}) become 
\begin{equation} 
	U_{t} + R\Lambda R^{-1}U_{x} = 0 \nonumber,
\end{equation}
\begin{equation} \label{euler5}
	R^{-1}U_{t} + \Lambda R^{-1}U_{x} = 0,
\end{equation}
where
\begin{equation} 
	R^{-1} = \left( \begin{array} {ccc} 
	1 & 0 & \frac{1}{K} \\ 
	0 & -\frac{1}{2 \sqrt{K}} & -\frac{1}{2K} \\ 
	0 & \frac{1}{2 \sqrt{K}} & -\frac{1}{2K} \end{array} \right), 
	\quad
	R = \left( \begin{array} {ccc} 
	1 & 1 & 1 \\ 
	0 & -\sqrt{K} & \sqrt{K}  \\ 
	0 & -K & -K \end{array} \right) \nonumber,
\end{equation}
\begin{equation} \label{euler6}
	\Lambda = V_{0}\left( \begin{array} {ccc} 
	0 &          &   \\ 
	  & \sqrt{K} &   \\ 
	  &          &  -\sqrt{K} \end{array} \right)	
\end{equation}
Based on the governing equations (\ref{euler5}) and (\ref{euler6}), we
have developed stable, particle-based, upwinding numerical schemes for
the system of Euler's equations. Details are described in the next
section.

\section{Numerical Discretization and Main Algorithms}

\subsection{Discrete Lagrangian Equtions} \label{section:DLE}

To solve numerically the hyperbolic system of PDE's (\ref{euler5}) and
(\ref{euler6}), the medium (compressible fluid or gas) is discretized
by a distribution of particles.  Each particle represents a Lagrangian
fluid cell of equal mass, and stores states of the continuum medium
such as density (that is proportional to the number density of
Lagrangian particles), pressure, internal energy, velocity, as well as
material properties and pointers to data structures containing material
models, such as the EOS.

To construct a Lagrangian upwinding scheme, we represent the system
(\ref{euler5})-(\ref{euler6}) in the following component-wise form
\begin{eqnarray} \label{sysEqn1}
	V_t + \frac{1}{K} P_t &=& 0,  \\
	-\frac{1}{2 \sqrt{K}} u_t - \frac{1}{2K} P_t
	 &=& -V_{0}\sqrt{K} \left[ -\frac{1}{2 \sqrt{K}} u_{x} - \frac{1}{2K} P_{x} \right] \label{sysEqn2}, \\
	\frac{1}{2 \sqrt{K}} u_t - \frac{1}{2K} P_t
	 &=& V_{0}\sqrt{K} \left[ \frac{1}{2 \sqrt{K}} u_{x} - \frac{1}{2K} P_{x} \right].  \label{sysEqn3}
\end{eqnarray} 
As $K>0$ for a thermodynamically consistent EOS, equation
(\ref{sysEqn2}) describes waves propagating from left to right, and
equation (\ref{sysEqn3}) describes waves propagating from right to
left. For an upwinding scheme, the spatial derivatives $u_x$ and $P_x$
will be computed on stencils within the corresponding physical domains of
dependence. Adding the subscripts $l$ and $r$ to the spatial derivatives
in equations (\ref{sysEqn2}) and (\ref{sysEqn3}), respectively, to
indicate that these terms, in the discrete form, will be computed using one-sided
derivatives, and solving for the temporal derivatives, we obtain
\begin{eqnarray}
	V_t &=& \frac{V_{0}}{2} \left( u_{xr}+u_{xl} \right) -
        \frac{V_{0}}{2 \sqrt{K}} \left( P_{xr}-P_{xl} \right), \label{eq15}\\ 
        u_t &=& \frac{V_{0} \sqrt{K}}{2} \left( u_{xr}-u_{xl} \right) -
        \frac{V_{0}}{2} \left( P_{xr}+P_{xl} \right), \label{eq16}\\
        P_t &=& -\frac{V_{0}K}{2} \left( u_{xr}+u_{xl} \right) + 
        \frac{V_{0} \sqrt{K}}{2} \left( P_{xr}-P_{xl} \right). \label{eq17}
\end{eqnarray}

An important component of a particle-based numerical scheme is the
calculation of differential operators based on states at the location
of particles. In Section 3.3, we describe in detail a method for both
numerical differentiation and interpolation based on local polynomial
fitting. In this section, we simply assume that we can compute
numerical approximations of differential operators with a desired
degree of accuracy on particle-based stencils located in the physical
domains of dependance.

The first-order ($O(\Delta t,\Delta x$) upwinding discretization of
the system (\ref{eq15}-\ref{eq16}) is obtained by the 1st order
discretization of spatial derivatives based on the local polynomial
fitting, and the 1st order discretization of temporal derivatives of
the state ($V,\, u$ or $P$) at the location of particle $j$,
\[
\frac{state_j^{n+1}-state_j^n}{\Delta t}.
\]
After the updates of states of each Lagrangian particle, particles are advanced by 
a mixture of the forward Euler scheme and backward Euler scheme:
\begin{eqnarray} \label{timeIntegrate}
	\frac{x^{n+1}-x^n}{\Delta t} &=& \frac{1}{2}\left(
        u^n+u^{n+1}\right)
\end{eqnarray}

The first order scheme is stable, provided that the standard CFL
condition is satisfied: $dt <= l/c$, 
where $l$ is the smallest interparticle
distance, but diffusive.  To reduce the amount of numerical diffusion of the 1st order scheme 
and obtain a higher order approximation on space and time,
we propose a \textit{modified Beam-Warming} scheme for the
Lagrangian particle system.  For the same reason as in the original
work on the Beam-Warming method \cite{BeamWarming}, an additional term
is added to equation (\ref{euler3}):
\begin{eqnarray}
	U_{t} + A(U)U_{x} - \frac{\Delta t}{2} A^2(U) U_{xx}= 0
        \nonumber \\ \Rightarrow U_t = -R \Lambda R^{-1} U_x +
        \frac{\Delta t}{2} R \Lambda^{2} R^{-1} U_{xx} \nonumber
        \\ \Rightarrow R^{-1} U_t = -\Lambda R^{-1} U_x + \frac{\Delta
          t}{2} \Lambda^{2} R^{-1} U_{xx}.
\end{eqnarray}
Equations (\ref{sysEqn2}) and (\ref{sysEqn3}) then become
\begin{eqnarray} \label{sysEqn6}
	 -\frac{1}{2 \sqrt{K}} u_t - \frac{1}{2K} P_t &=&
         -V_{0}\sqrt{K} \left[ -\frac{1}{2 \sqrt{K}} u_{xl} -
           \frac{1}{2K} P_{xl} \right] \nonumber \\ 
           &+& \frac{\Delta
           t}{2} V_{0}^2 K \left[ -\frac{1}{2 \sqrt{K}} u_{xxl} -
           \frac{1}{2K} P_{xxl} \right] ,\\ 
           \frac{1}{2 \sqrt{K}} u_t -
         \frac{1}{2K} P_t \label{sysEqn7} &=& V_{0}\sqrt{K} \left[
           \frac{1}{2 \sqrt{K}} u_{xr} - \frac{1}{2K} P_{xr} \right]
         \nonumber \\ 
         &+& \frac{\Delta t}{2} V_{0}^2 K \left[
           \frac{1}{2\sqrt{K}} u_{xxr} - \frac{1}{2K} P_{xxr} \right].
\end{eqnarray}
Solving equations (\ref{sysEqn1}), (\ref{sysEqn6}) and (\ref{sysEqn7}) yields
\begin{eqnarray}
	V_t &=& \frac{V_{0}}{2} \left( u_{xr}+u_{xl} \right) -
        \frac{V_{0}}{2 \sqrt{K}} \left( P_{xr}-P_{xl} \right),
        \nonumber \\ 
        &+& \frac{\Delta t}{4} \left[ V_{0}^2 \sqrt{K}
          \left( u_{xxr}-u_{xxl} \right) - V_{0}^2 \left(
          P_{xxr}+P_{xxl} \right)\right] ,\\ 
          u_t &=& \frac{V_{0}
          \sqrt{K}}{2} \left( u_{xr}-u_{xl} \right) - \frac{V_{0}}{2}
        \left( P_{xr}+P_{xl} \right), \nonumber \\ 
        &+& \frac{\Delta
          t}{4} \left[ V_{0}^2 K \left( u_{xxr}+u_{xxl} \right) -
          V_{0}^2 \sqrt{K} \left( P_{xxr}-P_{xxl} \right)\right], \\ 
        P_t &=& -\frac{V_{0}K}{2} \left( u_{xr}+u_{xl} \right) +
        \frac{V_{0} \sqrt{K}}{2} \left( P_{xr}-P_{xl} \right), \nonumber \\ 
        &+& \frac{\Delta t}{4} \left[ -V_{0}^2
          K^{\frac{3}{2}} \left( u_{xxr}-u_{xxl} \right) + V_{0}^2 K
          \left( P_{xxr}+P_{xxl} \right)\right].
\end{eqnarray}
By discretizing spatial derivatives using the second order local
polynomial fitting, as described in Section 3.3, we obtaine numerical
scheme that is second order in both time and space, $O(\Delta t^2,
\Delta x^2, \Delta t \Delta x)$, and conditionally stable. The CFL
condition is simlar to the one of the grid-based Beam-Warming scheme:
in 1D, $dt <= 2l/max(c,u)$. Note that time steps can be twice larger compared to the 
1st order scheme.
 
\subsection{Time Integration and Directional Splitting} \label{section:DS}
In this section, we focus on details of multidimensional schemes. We
present explicit formulas for equations in the three-dimensional space. The system
in the two-dimensional space is obtained by obvious reductions. 

In the three-dimensional space, the conservative form of the
Lagrangian formulation of the Euler equations is:
\begin{equation} \label{euler3d1}
U_t^{'} + \left[ F_1(U^{'}) \right]_x + \left[ F_2(U^{'}) \right]_y + \left[ F_3(U^{'}) \right]_z = 0,
\end{equation}
where
\begin{equation}
U^{'} = \left[ \begin{array} {c c c c c } V & u & v & w & E \end{array} \right]^T, \nonumber
\end{equation}
\begin{equation} \label{euler3d2}
F_1(U^{'}) = V_0 \left( \begin{array} {c} -u \\ P \\ 0 \\ 0 \\ Pu \end{array} \right), \, \,
F_2(U^{'}) = V_0 \left( \begin{array} {c} -v \\ 0 \\ P \\ 0 \\ Pv \end{array} \right), \, \,
F_3(U^{'}) = V_0 \left( \begin{array} {c} -w \\ 0 \\ 0 \\ P \\ Pw \end{array} \right).
\end{equation}
Assuming that the EOS is of the form $e=f\left( P,V \right) $ and
using $U=[V$ $u$ $v$ $w$ $P]^{T}$ as the state vector,
we can rewrite the equations in the following form
\begin{equation} \label{euler3d3}
	U_{t} + A_1 U_{x} + A_2 U_{y} + A_3 U_{z} = 0,	
\end{equation}
where
\begin{equation} \label{euler3d4}
	U = \left( \begin{array} {c} V \\ u \\ v \\ w \\ P\end{array} \right),
	\quad
	A_1 = V_{0}\left( 
	\begin{array} {ccccc} 
		0 & -1 & 0 & 0 & 0 \\ 
		0 & 0 & 0 & 0 & 1  \\
		0 & 0 & 0 & 0 & 0  \\
		0 & 0 & 0 & 0 & 0  \\ 
		0 & K & 0 & 0 & 0 
	\end{array} \right), \nonumber
\end{equation}
\begin{equation}	
	A_2 = V_{0}\left( 
	\begin{array} {ccccc} 
		0 & 0 & -1 & 0 & 0 \\ 
		0 & 0 & 0 & 0 & 0  \\
		0 & 0 & 0 & 0 & 1  \\
		0 & 0 & 0 & 0 & 0  \\ 
		0 & 0 & K & 0 & 0 
	\end{array} \right),
	\quad
	A_3 = V_{0}\left( 
	\begin{array} {ccccc} 
		0 & 0 & 0 & -1 & 0 \\ 
		0 & 0 & 0 & 0 & 0  \\
		0 & 0 & 0 & 0 & 0  \\
		0 & 0 & 0 & 0 & 1  \\ 
		0 & 0 & 0 & K & 0 
	\end{array} \right),
\end{equation}
where $K$ is defined in equation (\ref{Kdef}).
We solve the system of hyperbolic PDEs (\ref{euler3d3} - \ref{euler3d4}) by
using the directional splitting method by Strang \cite{Strang_splitting}.  
Specifically, instead of solving equation (\ref{euler3d3}),
one solves separately the following three system of PDEs:
\begin{eqnarray} \label{euler3d5}
	U_t + 3A_1 U_x = 0,\\  \label{euler3d6}
	U_t + 3A_2 U_y = 0,\\ \label{euler3d7}
	U_t + 3A_3 U_z = 0,
\end{eqnarray} 
which is equivalent to solving 
\begin{eqnarray} \label{euler3d8}
	U_{1t} + 3A U_{1x} = 0, \\ \label{euler3d9}
	U_{2t} + 3A U_{2y} = 0, \\ \label{euler3d10}
	U_{3t} + 3A U_{3z} = 0,
\end{eqnarray}
where
\begin{equation} \label{euler3d11}
	U_1 = \left( \begin{array} {c} V \\ u \\ P \end{array} \right), 
	\,
	U_2 = \left( \begin{array} {c} V \\ v \\ P \end{array} \right),
	\,
	U_3 = \left( \begin{array} {c} V \\ w \\ P \end{array} \right),
	\,
	A = V_0 \left( \begin{array} {ccc} 
	0 & -1 & 0 \\
	0 &  0 & 1 \\
	0 &  K & 0
	\end{array} \right).
\end{equation}
Each of the three system of equations 
(\ref{euler3d8}) - (\ref{euler3d10})
is solved by the techniques introduced 
in section \ref{section:DLE}, and
the solutions are combined in the following order
\begin{equation}
	\left( \begin{array} {c} V \\ u \\ v \\ w \\ p \end{array} \right)^{(0)}
	\underrightarrow{(\ref{euler3d8})}
	\left( \begin{array} {c} V \\ u \\ v \\ w \\ p \end{array} \right)^{(\frac{\Delta t}{6})}
	\underrightarrow{(\ref{euler3d9})}
	\left( \begin{array} {c} V \\ u \\ v \\ w \\ p \end{array} \right)^{(\frac{2\Delta t}{6})}
	\underrightarrow{(\ref{euler3d10})}
	\left( \begin{array} {c} V \\ u \\ v \\ w \\ p \end{array} \right)^{(\frac{4\Delta t}{6})} \nonumber	
\end{equation}
\begin{equation} \label{SSplit3d}
	\underrightarrow{(\ref{euler3d9})}
	\left( \begin{array} {c} V \\ u \\ v \\ w \\ p \end{array} \right)^{(\frac{5\Delta t}{6})}
	\underrightarrow{(\ref{euler3d8})}
	\left( \begin{array} {c} V \\ u \\ v \\ w \\ p \end{array} \right)^{(\Delta t)}
\end{equation}
where $\Delta t$ denotes one discrete time step satisfying the CFL condition. 
The Strang splitting method maintains the second order of accuracy if the accuracy of each step is not lower than second,
making it unnecessary for the 1st order numerical scheme.
To implement the modified Beam-Warming scheme within the Strang splitting
 steps (\ref{euler3d8}) - (\ref{euler3d10}), we solve  the following equations
\begin{eqnarray} \label{euler3dBW4}
	U_{1t} + 3\left(A U_{1x} -\frac{\Delta t}{2} A^2 U_{1xx} \right)= 0, \\  \label{euler3dBW5}
	U_{2t} + 3\left(A U_{2y} -\frac{\Delta t}{2} A^2 U_{2yy} \right)= 0, \\  \label{euler3dBW6}
	U_{3t} + 3\left(A U_{3z} -\frac{\Delta t}{2} A^2 U_{3zz} \right)= 0.
\end{eqnarray}
The solutions to equations (\ref{euler3dBW4}) - (\ref{euler3dBW6}) are then combined 
by equation (\ref{SSplit3d}) to obtain the complete second order solution to equation
(\ref{euler3d3}).

\subsection{Local Polynomial Fitting}
\label{section:LPF}

The local polynomial fitting on arbitrary sets of points has long been
used to obtain approximation of functions and their derivatives.
Details of the method and its accurracy is discussed in
\cite{BenitoUrena01,Fan93localpolynomial,Jiao08}.  Generally, 
$\nu$th order derivative can be approximated  with $(n-\nu+1)$th order of accuracy 
using $n$th order polynomial. For simplicity, a 2D example is discussed here.  In the vicinity of
a point $0$, the function value in the location of a point $i$ can be expressed by
the Taylor series as
\begin{equation}
U_i = U_0+h_i\left.\frac{\partial U}{\partial x}\right|_0+k_i\left.\frac{\partial U}{\partial y}\right|_0+\frac12\left(h^2_i\left.\frac{\partial^2 U}{\partial x^2}\right|_0 + k^2_i\left.\frac{\partial^2 U}{\partial y^2}\right|_0 + 2h_ik_i\left.\frac{\partial^2 U}{\partial x \partial y}\right|_0\right) + \ldots ,
\end{equation}
where $U_i$ and $U_0$ are the corresponding function values in the location of points
$i$ and $0$,  $h_i = x_i - x_0$, $k_i = y_i - y_0$, 
and the derivatives are calculated in the location of the point $0$.
A polynomial can be used to approximate the original
function and we employ a second order polynomial in this example:
\begin{equation}
\tilde{U}=U_0+h_i\theta_1+k_i\theta_2+\frac12h^2_i\theta_3+\frac12k^2_i\theta_4+h_ik_i\theta_5 .
\end{equation} 
Here the variables $\theta_1$, $\theta_2$, $\theta_3$, $\theta_4$ and
$\theta_5$ are the estimates for $ \frac{\partial U}{\partial x}$,
$\frac{\partial U}{\partial y}$, $\frac{\partial^2 U}{\partial x^2}$,
$\frac{\partial^2 U}{\partial y^2}$, and $\frac{\partial^2 U}{\partial
  x \partial y}$, respectively.  In order to compute values of these variables,
we perform a local polynomial fitting using $m>=5$ points in the
vicinity of center point $0$.  The following linear system $Ax=b$
\begin{equation} \label{linear_system_axb}
\begin{bmatrix} h_1 & k_1 & \frac12 h^2_1 & \frac12 k^2_1 & h_1k_1\\ 
                h_2 & k_2 & \frac12 h^2_2 & \frac12 k^2_2 & h_2k_2 \\  
                \vdots & \vdots & \vdots & \vdots & \vdots \\ 
                h_n & k_n & \frac12 h^2_n & \frac12 k^2_n & h_nk_n
\end{bmatrix}  
\left[ \begin{array}{c} 
\theta_1 \\  \theta_2 \\ \theta_3 \\ \theta_4 \\ \theta_5 
\end{array} \right] 
=
\left[ \begin{array}{c} 
U_1 - U_0 \\ U_2 - U_0 \\ \vdots \\ U_n - U_0 
\end{array} \right],
\end{equation}
is usually overdetermined.
As a proper selection of a neighborhood is important for accuracy and stability,
neighbor search algorithms used in our upwind solvers are described in the next
subsection.

An optimal solution to (\ref{linear_system_axb}) is a solution $x$ that minimizes the $L_2$ norm of
the residual, i.e.,
\begin{equation}
min\|Ax-b\|_2,
\end{equation}
and the QR decomposition with column pivoting is employed to obtain $x$. Suppose
\begin{equation}
A=Q\begin{bmatrix} R\\0\end{bmatrix}P^T, m\ge n,
\end{equation}
where $Q$ is an orthonomal matrix, $R$ is an upper triangle matrix, 
and $P$ is a permutation matrix, chosen (in general) so that
\begin{equation}
|r_{11}|\ge|r_{22}|\ge\cdots\ge|r_{nn}|.
\end{equation}
Moreover, for each $k$,
\begin{equation}
|r_{kk}|\ge\|R_{k:j,j}\|_2 
\end{equation}
for $j = k+1,\cdots,n$. One can numerically determine an index k, 
such that the leading submatrix $R_{11}$ in the first $k$ rows and columns 
is well conditioned and $R_{22}$ is negligible:
\begin{equation}
R = \begin{bmatrix}R_{11}&R_{12}\\0&R_{22}\end{bmatrix}\simeq\begin{bmatrix}R_{11}&R_{12}\\0&0\end{bmatrix}
\end{equation}
Then $k$ is the effective rank of A. 
Discussion about the numerical rank determination can be found in \cite{GolubVanLoan96}. 
A simple way to determine numerical rank is to set a tolerance $\epsilon$ 
and find the first $k$ such that
\begin{equation} \label{effective_rank}
R_{kk} < \epsilon R_{11}.
\end{equation}
If there is such $k$, then the effective numerical rank is $k-1$.  The
choice of $\epsilon$ is $10^{-3}$ in many of the test problems
discussed in later sections.  The solution for linear system
(\ref{linear_system_axb}) can be obtained as
\begin{equation}
x=P\left[ \begin{array}{c}R_{11}^{-1}c_1 \\ 0 \end{array}\right]
\end{equation}
where $c_1$ is the first $k$ elements of $c=Q^Tb$. This can also be written as
\begin{equation}
x=A^+b,
\end{equation}
where
\begin{equation}
A^+=P\begin{bmatrix}R_{11}^{-1}&0\\0&0\end{bmatrix}Q^T
\end{equation}
is the pseudoinverse of matrix $A$.

\subsection{The neighbor Search Algorithm and Dynamic Stencil Selection}
In a simulation involving $N$ Lagrangian particles, a new
\textit{stencil} of neighbors - those particles used for solving
equation (\ref{linear_system_axb}) - is to be selected at the begining
of each time step for all the $N$ particles.  As a result, it is
critical that an efficient neighbor search algorithm is employed. 
The neighbor search method is described in the next subsection.
While the neighbor search is to obtain a group of particles lying
within some pre-specified distance away from the particle of interest,
it is not necessary that all these neighbors are used in numerical stencils.
The selection of stencil points from neighbors to ensure the accuracy and stability
is discussed in Subsection 3.4.2.

\subsubsection{Neighbor Search Algorithms}

One of the main advantagees of the Lagrangian particle method compared to grid-based methods
is its ability to simulate large and extremely non-uniform domains. By a non-uniform domain we mean
a domain in which only a small fraction of the total volume
occupied by matter, found typically in astrophysics and high energy density physics, and other applications
dealing with dispersed fragments of matter.
For these applications, we  use $2^k$-tree neighbor search algorithms \cite{Samet_octree}.
The $2^k$-tree is a tree data structure in a k-dimensional space in which each node has at most
$2^k$ dependents.  Quadtree and octree  are the standard terms in 2D and 3D spaces, respectively.
The tree construction can be performed with $O(N\log N)$ operation.  In this process, the choice of the tree depth is
essential and the optimal empirical number is four or five.  After the construction step,
the search of a tree for obtaining the neighborhood of a particle can be performed with $O(\log N)$ operation.

However the $2^k$-tree method is not universally optimal for all types of problems.
If the computational domain is almost uniformly filled with a weakly compressible matter in which inter-particle distances
change insignificantly during the simulation allowing the use of the same neighbor search radius for all particles, the search
of neighbors can be performed in constant time.  In this case, we use the bucket search algorithm \cite{Monaghan85}.  
The entire computational domain is divided into square (cubic) cells of the side length equal to the search radius $r$.
 For each particle inside a cell, only the neighboring cells need to be considered in the search process. Clearly, the method is
 not optimal if the location of matter in the space is very sparse, and is not applicable if different search radii must be used for 
 different particles. The $2^k$-tree neighbor search algorithm is more universal and applicable to a wide range
of problems. In a forthcoming paper, we will describe an optimal parallel octree neighbor search algorithm for a 3D Lagrangian particle code. 

\subsubsection{Dynamic stencil selection}

After the neighbor search step, each particle obtains a list of neighbors
which lie within the range of a pre-specified search radius.  
To enforce upwinding, however, only one-sided information should be
used when solving equation (\ref{linear_system_axb}).  
For the calculation of one-sided derivatives, each particle will, in general, have four one-sided
neighborhoods in two-dimensions, and six neighborhoods in three-dimensions.

Without loss of generality, the process of the dynamic stencil selection will
be discussed using an example of computing $u_{xr}$.
After gathering one-sided neighbors, two main issues must be resolved for
accurate evaluation of spatial derivatives.  The first one
is related to the \textit{shape} of the neighborhood.  The list of
one-sided neighbors is sorted by their distance from the center
particle in ascending order.  Suppose we obtain the following sorted list of
neighbors of the particle $0$ for computing $u_{xr}$:
\[
	\left\{ p_{1u}, \, p_{2u}, \, p_{3l}, \, p_{4u}, \, p_{5u}, 
	       \, p_{6u}, \, p_{7l}, \, p_{8u}, \, p_{9l} \right\} .
\]
Here the subscripts $u$ and $l$ indicate the upper and lower half-planes in the $y$-direction:
$y_i >= y_0$ and $y_i < y_0$, respectively.  If a simple distance-based algorithm picks up six neighbors,
then the corresponding stencil is composed of
\begin{equation}
	\left\{ p_{1u}, \, p_{2u}, \, p_{3l}, \, p_{4u}, \, p_{5u}, \,
        p_{6u} \right\}, \nonumber
\end{equation}
thus producing a highly \textit{unbalanced} stencil in terms of the \textit{shape}. 
Therefore, besides sorting neighbors in ascending order of the
distance from the center particle, the order is rearranged such that
neighbors from the upper half and lower half occur interchangeably in
the list
\begin{equation}
	\left\{ p_{1u}, \, p_{3l}, \, p_{2u}, \, p_{7l}, \, p_{4u}, 
	       \, p_{9l}, \, p_{5u}, \, p_{6u}, \, p_{8u} \right\} \nonumber
\end{equation}
The six-particle stencil now becomes:
\begin{equation}
	\left\{ p_{1u}, \, p_{3l}, \, p_{2u}, \, p_{7l}, \, p_{4u}, \,
        p_{9d} \right\} \nonumber
\end{equation}
This approach yields more balanced-in-shape stencils, 
and typically results in more accurate spatial derivatives.

The second issue is the optimization of the $\it number$ of neighbors 
for solving equation (\ref{linear_system_axb}).
In the case of second order local polynomial fitting, for example,
five neighbors are required to solve equation (\ref{linear_system_axb}).
However, as equation (\ref{effective_rank}) suggests, 
the effective rank of the matrix $A$
in equation (\ref{linear_system_axb}) may be less than five if
\begin{equation} 
 R_{kk} < \epsilon  R_{11}, \, k=2,3,4,5.
\end{equation}
To avoid rank deficiency, a dynamic process for selecting neighbors
into the stencil is designed.  First, select the \textit{tolerance}
parameter $\epsilon$ as in equation (\ref{effective_rank}).  For the
case of second order local polynomial fitting, one starts with five or
six neighbors in the stencil. Based on this stencil, the QR decomposition
with column pivoting is performed.  Then determine the effective rank
by equation (\ref{effective_rank}).  If the effective rank is no less
than five, the stencil is complete. Otherwise, the next neighbor in the
neighbor list is added to the stencil. The process continues until  the
effective rank is regained. For instance, in some cases one may need
to use seven neighbors to gain an effective rank of five. 
Our algorithms uses much smaller number of neighbor particles
compared to the modified versions of SPH such as Godunov-SPH,
P-SPH, PHANTOM etc. that may require hundreds of particles \cite{GIZMO}.

In certain cases, the neighbor list may not contain sufficient number of 
particles needed by a stencil. In such a case,
one may consider lowering the order of local polynomial fitting for
this particle in the given direction.  Lowering to first order local polynomial
fitting requires only an effective rank of two for the two-dimensional
case. 

\subsection{Limiters}
\label{sec:LIM}

The second order Lagrangian particle algorithm based on the modified 
Beam-Warming scheme is dispersive. 
To eliminate the resulting oscillations, a new type of 
limiter based on divided difference was developed and  coupled with the numerical 
integration. The application of the algorithm with the limiter is
demonstrated in section \ref{section:num_res}.

In the flux-limiter method \cite{} TODO, the magnitude of the correction depends on 
the smoothness of data (represented by $\Phi$), and can be written as
\begin{equation}
	F(U;j) = F_L(U;j) + \Phi(U;j) [F_H(U;j)-F_L(U;j)]
	\label{eq:lim1}
\end{equation}
In order to measure the smoothness of data,
we can use the ratio of consecutive gradients:
\begin{equation}
	\theta_j = \frac{U_j-U_{j-1}}{U_{j+1}-U_j}
	\label{eq:lim2}
\end{equation}
Or we can use the average of the ratio of consecutive gradients from both directions:
\begin{equation}
	\theta_j = \frac{1}{2} \left( \frac{U_j-U_{j-1}}{U_{j+1}-U_j} + \frac{U_{j+1}-U_j}{U_j-U_{j-1}} \right)
	\label{eq:lim2prime}
\end{equation}
Equation \eqref{eq:lim2prime} has the advantage that it is symmetric.
If $\theta_j$ is near $1$ the data is presumably smooth.
If $\theta_j$ is far from $1$ there may be discontinuity
near data $U_j$.
Let $\Phi(U;j)\equiv\phi_j$ to be a function of $\theta_j$:
\begin{equation}
	\phi_j=\phi(\theta_j)
	\label{eq:lim3}
\end{equation}
Van Leer \cite{VanLeer74} proposed a smooth limiter function
\begin{equation}
	\phi(\theta)=\frac{|\theta|+\theta}{1+|\theta|}
	\label{eq:lim4}
\end{equation}
We let $\phi(\theta)=0$ for $\theta < 0$
or when $\theta$ is arbitrarily large.
Nore that $\theta < 0$ in the case when
$U_{j+1}-U_j$ and $U_j-U_{j-1}$ are in opposite signs
in both equations \eqref{eq:lim2} and \eqref{eq:lim2prime}).
$\theta$ is arbitrarily large when 
$U_{j+1}-U_j=0$ in equation \eqref{eq:lim2}
and when $U_{j+1}-U_j=0$ or $U_j-U_{j-1}=0$ 
in equation \eqref{eq:lim2prime}.

Without loss of generality, we demonstrate the idea using 
the flux for volume. Remind that in the proposed Lagrangian 
particle method the volume flux is defined as 
(equation \eqref{eq:flux1}) 
\begin{align}
	V_t &= \frac{V_{0}}{2} \left( u_{xr}+u_{xl} \right) -
        \frac{V_{0}}{2 \sqrt{K}} \left( P_{xr}-P_{xl} \right)
        \nonumber \\ &+ \frac{\Delta t}{4} \left[ V_{0}^2 \sqrt{K}
          \left( u_{xxr}-u_{xxl} \right) - V_{0}^2 \left(
          P_{xxr}+P_{xxl} \right)\right] 
	\label{eq:flux1_} 
\end{align}
Let the lower order flux (first order flux) of volume be defined as
\begin{equation} 
	F_L(V) = \frac{V_{0}}{2} \left( u_{xr(1)}+u_{xl(1)} \right) -
        \frac{V_{0}}{2 \sqrt{K}} \left( P_{xr(1)}-P_{xl(1)} \right)
		\label{eq:LF}
\end{equation}
where the subscript $(1)$ denotes the spatial derivatives obtained
by first order polynomial fitting.
Then define the higher order flux (second order flux) of volume be defined as
\begin{align} 
	F_H(V) &= \frac{V_{0}}{2} \left( u_{xr(2)}+u_{xl(2)} \right) -
        \frac{V_{0}}{2 \sqrt{K}} \left( P_{xr(2)}-P_{xl(2)} \right)
        \nonumber \\ &+ \frac{\Delta t}{4} \left[ V_{0}^2 \sqrt{K}
          \left( u_{xxr(2)}-u_{xxl(2)} \right) - V_{0}^2 \left(
          P_{xxr(2)}+P_{xxl(2)} \right)\right] 
	\label{eq:HF} 
\end{align}
where the subscript $(2)$ denotes the spatial derivatives obtained
by second order polynomial fitting.

In order to make the measure of the smoothness of data ($\theta$)
generalizable to higher dimensions and applicable to the 
Lagrangian particle mehtod,
we propose using the one-sided spatial derivatives calculated 
by methods introduced in section \ref{section:LPF}.
Depneding on the type of data we have, $\theta$ of a particle $j$
is calculated as:
\begin{align}
	\theta_j(u) &= \frac{u_{xl(1)}}{u_{xr(1)}} \label{eq:theta_u1} \\
	\theta_j(P) &= \frac{P_{xl(1)}}{P_{xr(1)}} \label{eq:theta_P1}
\end{align}
where $u$ is the velocity in the $x$, $y$, or $z$-direction and $P$ is the pressure.
Alternatively, we can also use
\begin{align}
	\theta_j(u) &= \frac{1}{2} \left( \frac{u_{xl(1)}}{u_{xr(1)}} + \frac{u_{xr(1)}}{u_{xl(1)}} \right) 
	\label{eq:theta_u2} \\
	\theta_j(P) &= \frac{1}{2} \left( \frac{P_{xl(1)}}{P_{xr(1)}} +  \frac{P_{xr(1)}}{P_{xl(1)}} \right) 
	\label{eq:theta_P2}
\end{align}
Note that equations \eqref{eq:theta_u2} and \eqref{eq:theta_P2} 
are better over \eqref{eq:theta_u1} and \eqref{eq:theta_P1} since they are symmetric. 
We choose 
\begin{equation}
	\theta_j = \min(\theta_j(u),\theta_j(P))
	\label{eq:thetaval}
\end{equation}
and calculate $\phi_j$ by \eqref{eq:lim4}. 
Note that we set $\phi_j=0$ 
when $\theta_j < 0$ or $\theta_j$ is arbitrarily large.
Substituting the calculated $\phi_j$, the lower and higher
flux in \eqref{eq:LF} and \eqref{eq:HF} into \eqref{eq:lim1},
we obtain the volume flux of particle $j$ as:
\begin{equation}
	F_j(V) = F_L(V;j) + \phi_j [F_H(V;j)-F_L(V;j)]
\end{equation}
Then time integration gives the volume at next time step as:
\begin{equation}
	V^{n+1} = V^{n} + \Delta t F_j(V)
\end{equation}

Similarly to  \cite{Shu_WENO}, we use the divided differences
to detect the region that contains discontinuities.
However, while the method of divided differences is used in \cite{Shu_WENO} for choosing 
one of different stencils, it is employed in our work for switching between 
high order and low order discretization schemes. 

The proposed limiter works as a  
switch between higher and lower order schemes to avoid the oscillatory behavior.

\subsection{Modelling of Free Surfaces using Ghost Particles}
\label{section:ghost_particle}
An important feature of the Lagrangian particle method is its ability
to robustly handle free surface flows with geometrically complex interfaces. 
The method is also generalizable to multiphase problems.
Here by free surface flows we mean  flows of fluid or gas in vacuum, and by multiphase problem we
mean the interface dynamics between two immiscible fluids or gases.  
In this section, we describe an algorithm for physically consistent solutions at free fluid or gas interfaces.  

The fluid / vacuum interface is modeled in our method by using \textit{ghost} particles in the vacuum region. 
A geometric algorithm places patches of ghost particles outside the fluid boundary, ensures their proper 
distance to the interface, and eliminates those particles that were placed too closely or inside the fluid.
Then the ghost particles are assigned physics states. 
The only functionality of ghost particles is to serve as neighbors of fluid particles when calculating spatial derivatives.
Hence, only two states are relevant: pressure and velocity.
As ghost particles represent vacuum, their pressure state is assigned to zero.
A weighted $0$th order local polynomial fitting is used to assign velocity states to ghost particles.
This involves computing the weighted
average of velocities of the fluid particles that are in a neighborhood of the ghost particle.
Let's assume that the weighting function of the particle $0$
in a three-dimensional space is $w(h_j,k_j,g_j)$, where $h_j=x_j-x_0$, $k_j=y_j-y_0$,
$g_j=z_j-z_0$, and $j$ is the index of neighbor particles. The velocity $u_0$ of particle $0$ satisfies
\begin{equation}
	\min_{u_0}\sum_{j=1}^N\left[ (u_0 - u_j)w(h_j,k_j,g_j)\right]^2,
\end{equation}
which leads to the solution 
\begin{equation}
	u_0 = \frac{\sum_{j=1}^N u_j w_j^2}{\sum_{j=1}^N w_j^2}
\end{equation} 
This simple algorithm adequately handles  the fluid / vacuum interface, but a Riemann solver-based algorithm will be used for  interfaces in
multiphase problems.

\section{Numerical Results}
\label{section:num_res}
In this section, we present results of one- and two-dimensional simulation
that serve as verification tests for the Lagrangian particle method, including the
free surface algorithm.

\subsection{1D Gaussian Pressure Wave Propagation with Periodic Boundaries}

\begin{figure}
	\centering
	\includegraphics[scale=0.9]{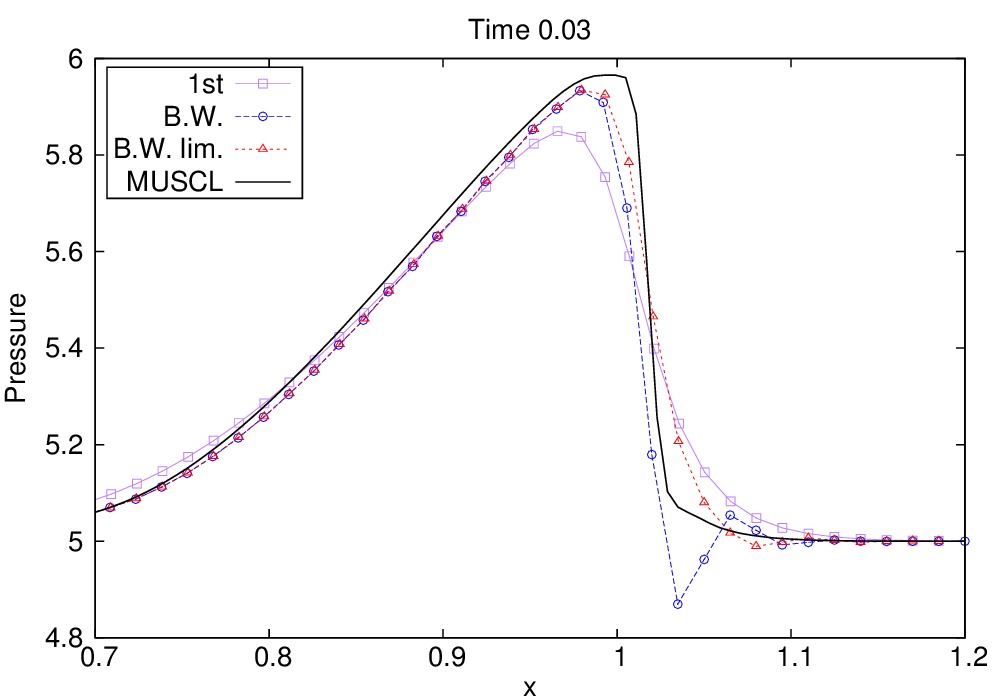} 
	\includegraphics[scale=0.9]{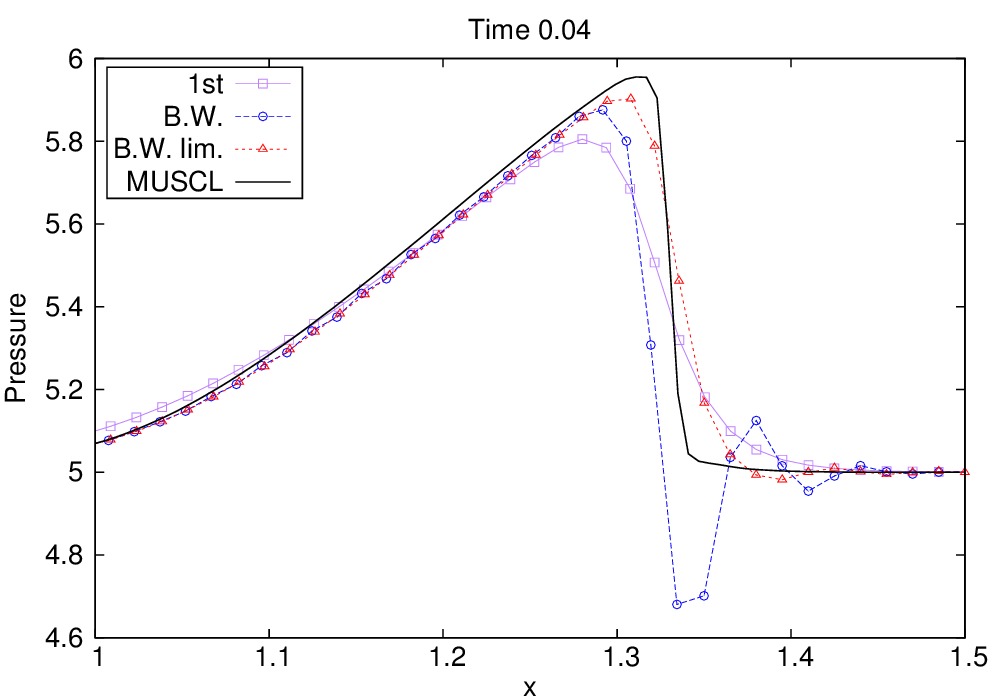} 
	\caption{Gaussian pressure wave propagation with periodic boundaries at time 0.03 (top) and 0.04 (bottom).
	Coarse-resolution simulations results were used to illustrate the behavior qualitatively.}
	\label{fig:gauss_poly_period_1d_t003_}
\end{figure}

We study the propagation of a pressure wave in gas with the constant initial density $\rho = 0.01$ and
the initial Gaussian pressure distribution
\begin{equation}
	p=5+2e^{-100x^2}  
\end{equation}
in the domain $-1.5\leq x\leq1.5$ with periodic boundaries on both ends.  The  polytropic gas EOS is 
used with $\gamma=5/3$.    The goal of the simulation is to demonstrate the accuracy 
of the proposed algorithm in resolving nonlinear waves with the formation of shocks. 
The benchmark data is obtained using a highly refined, grid-based  MUSCL scheme.
The results, shown in Figure \ref{fig:gauss_poly_period_1d_t003_}, are labeled as \textit{1st} for the first order
local polynomial fitting, \textit{B.W.} for the
Beam-Warming scheme with second order local polynomial fitting, and \textit{B.W. lim.} for the Beam-Warming
scheme with the second order local polynomial fitting with limiter, respectively.
As expected, first order scheme is diffusive,
while the Beam-Warming scheme is dispersive near discontinuities.
However, results demonstrates that the proposed limiter method
effectively reduces dispersions near sharp edges, resulting in maintaining globally the second order of convergence.

We have also verified that the Lagrangian particle methods accurately resolves waves in stiff materials.  We use the same
initial pressure as before, but replace the polytropic EOS with the stiffened polytropic EOS
\begin{equation} \label{spolyEOS}
	E = \frac{\left( P+\gamma P_{\infty}\right)V}{\gamma -1}
\end{equation}
with $\gamma = 6$ and $P_{\infty} = 7000$.
The convergence results can be found in Tables {\ref{tab:ConvPoly}} and 
{\ref{tab:ConvSpolyCFL0.5}}. In both cases, the second order convergence is obtained.
\begin{table}
\begin{center}
    \begin{tabular}{| l | l | l |}
    \hline
    Number of particles & Relative $L2$-norm error & Rate of Convergence\\ \hline
    $240$  & $0.051$   &  NA         \\ \hline
    $480$  & $0.018 $   &  $2.88$  \\ \hline
    $960$  & $0.0049  $  &  $3.60$ \\ \hline
    $1920$ & $0.0012  $   &  $4.02$ \\ \hline
    $3840$ & $0.00029  $ &  $4.23$ \\ \hline
    $7680$ & $0.000068$ &  $4.24$ \\ \hline
    \end{tabular}\\
\end{center}
\caption{Convergence for the polytropic gas EOS case with 
$\gamma=\frac{5}{3}$ and initial density 
$\rho_0 = \frac{1}{V_0} = 0.01$}
\label{tab:ConvPoly}
\end{table}
\begin{table}
\begin{center}
    \begin{tabular}{| l | l | l |}
    \hline
    Number of particles & Relative $L2$-norm error & Rate of convergence\\ \hline
    $240$  & $0.069$      &  NA         \\ \hline
    $480$  & $0.021$       &  $3.23$ \\ \hline
    $960$  & $0.0056$     &  $3.84$ \\ \hline
    $1920$ & $0.0014$     &  $3.97$ \\ \hline
    $3840$ & $0.00035$    &  $4.0$ \\ \hline
    $7680$ & $0.000093$    &  $3.8$ \\ \hline
    \end{tabular}\\
\end{center}
\caption{Convergence for the stiffened polytropic gas EOS case with 
$\gamma=6$, $P_{\infty}=7000$, and initial density 
$\rho_0 = \frac{1}{V_0} = 1$}
\label{tab:ConvSpolyCFL0.5}
\end{table}

\subsection{2D Gaussian Pressure Wave Propagation with Free Surface}
\begin{figure}
	\centering
	\includegraphics[width=0.32\textwidth]{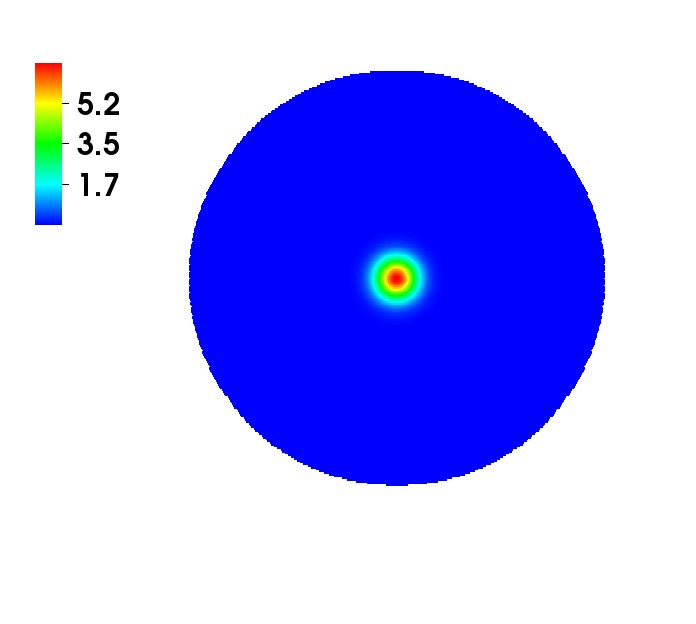} 
	\includegraphics[width=0.32\textwidth]{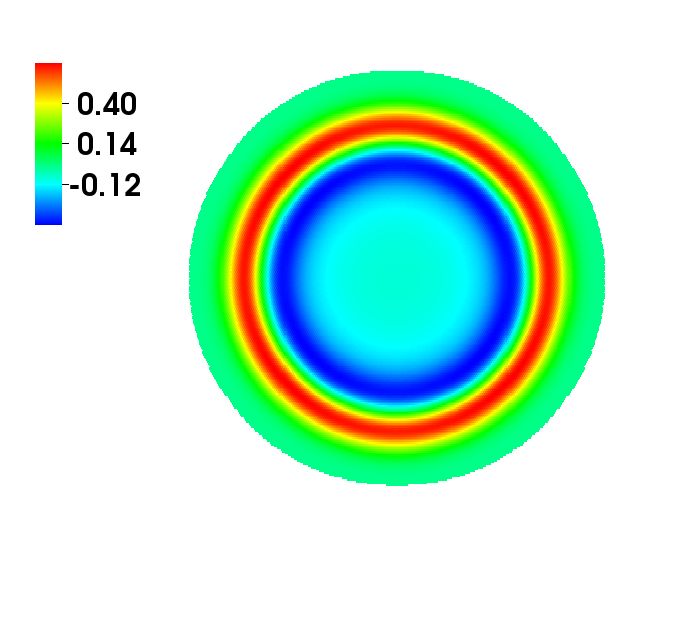} 
	\includegraphics[width=0.32\textwidth]{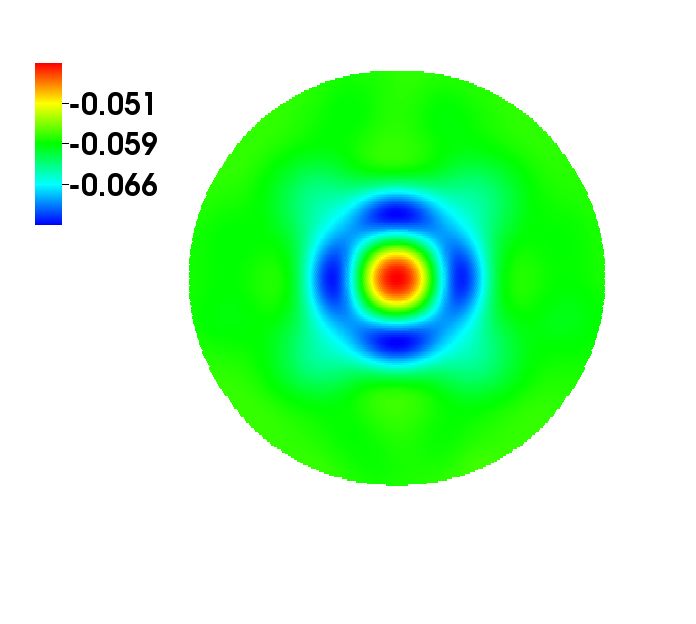} 
	\caption{2D Gaussian pressure wave propagation in disk with free surface. Pressure distribution (bar)
at initial time (left), 10 (middle), and 60  (right).}
	\label{fig:gauss_spoly_free_2d_0}
\end{figure}

To  test the proposed algorithm for
two-dimensional problems involving free surfaces, a circular disk of particles
with stiffened polytropic gas EOS (with $\gamma=6$,
$P_{\infty}=7000$, and $\rho=1$) and a Gaussian
pressure profile was initialized.  The results are presented in two
dimensions in Figure \ref{fig:gauss_spoly_free_2d_0}. Note that the
latest-time plot in Figure \ref{fig:gauss_spoly_free_2d_0} represents
the state when the pressure waves have been reflected from the oscillatory free
surface for more than ten times. To verify the
accuracy, the analogous one-dimensional problem with cylindrical
coordinates under the Eulerian formulation was solved using a refined MUSCL
scheme with the method of front tracking for the free surface implemented in the 
FronTier code \cite{FTcode}. The location and shape of the pressure wave and the interface
as well as the oscillatory motion of the free surface are
 in good agreement with the FronTier simulation.
The verification test and  the fact that the
pressure wave maintains good
symmetry after many reflections from the free surface demonstrate that 
the method for modeling vacuum introduced in section \ref{section:ghost_particle}
works well with the proposed algorithm.  


\subsection{2D Shock Tube Problem}

\subsection{Gresho Vortex}

The Gresho problem is a steady-state rotating vortex roblem that has an exact analytic solution 
in the case of Euler equations \cite{LiskaWen03}. Gresho vortex is a, inviscid gas vortex with 
such a radial distribution of the angular velocity $u_{\phi}$ and pressure, that the centripital
force is compensated by the gradient of pressure:
\begin{equation}
\left( u_{\phi}(r),p(r)\right) = \left\{
\begin{array}{ll}
\left( 5r, 5+\frac{25}2 r^2 \right), & 0\le r\le 0.2,\\
\left( 2-5r, 9-4\ln0.2 + \frac{25}2r^2 - 20r + 4\ln r \right), & 0.3\le r\le 0.4,\\
\left( 0, 3+4\ln 2 \right), & 0.4\le r.
\end{array}
\right.
\end{equation}
 
\begin{figure}
\centering
\subfigure[Pressure, t=0]{\label{gresho_lp_p0}\includegraphics[width=0.49\textwidth]{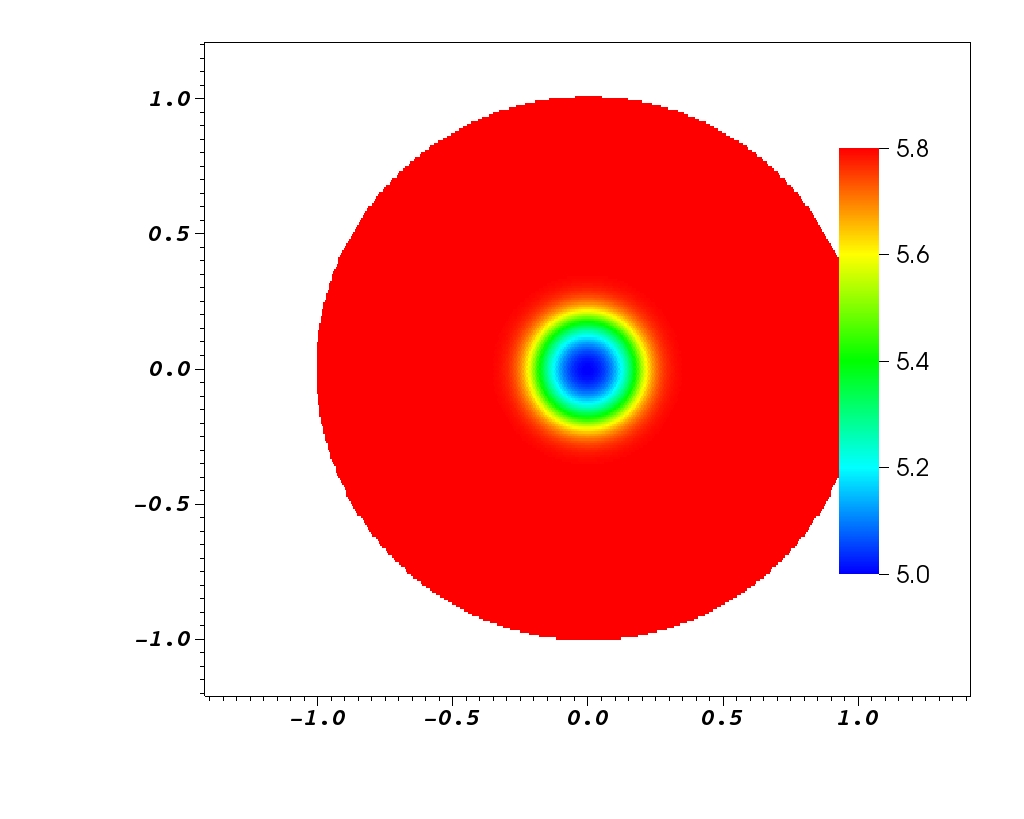}}
\subfigure[Pressure, t=1]{\label{gresho_lp_p1}\includegraphics[width=0.49\textwidth]{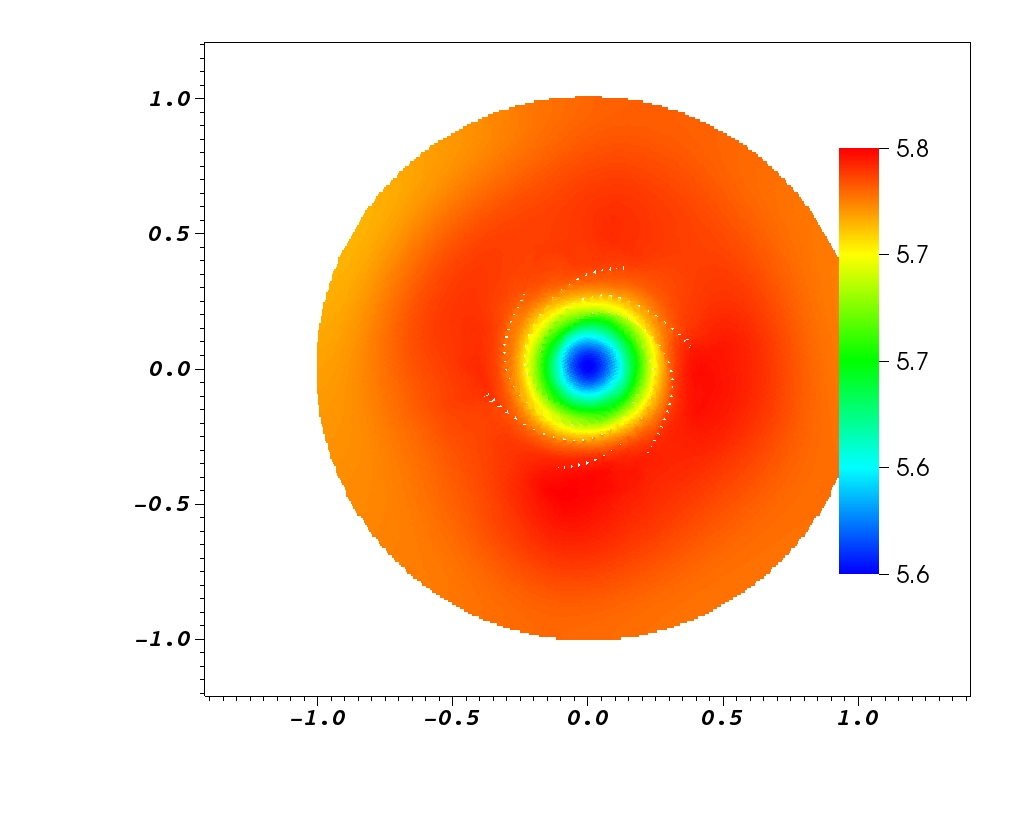}}
\subfigure[Velocity, t=0]{\label{gresho_lp_u0}\includegraphics[width=0.49\textwidth]{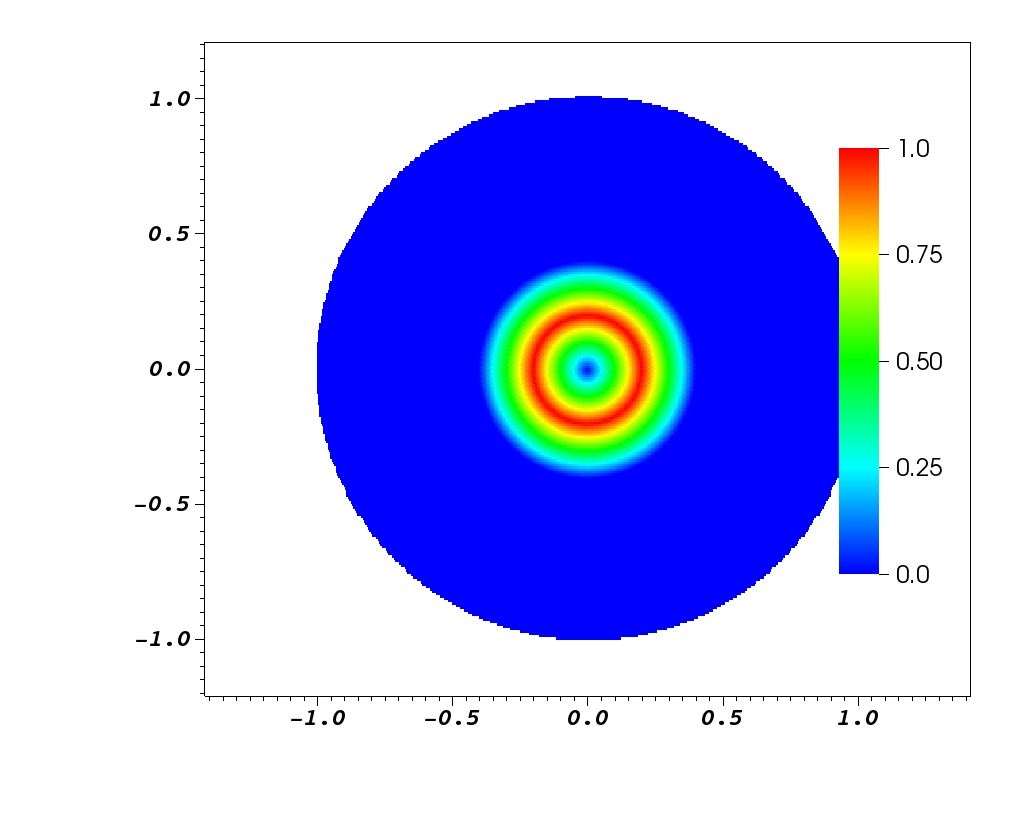}}
\subfigure[Velocity, t=1]{\label{gresho_lp_u1}\includegraphics[width=0.49\textwidth]{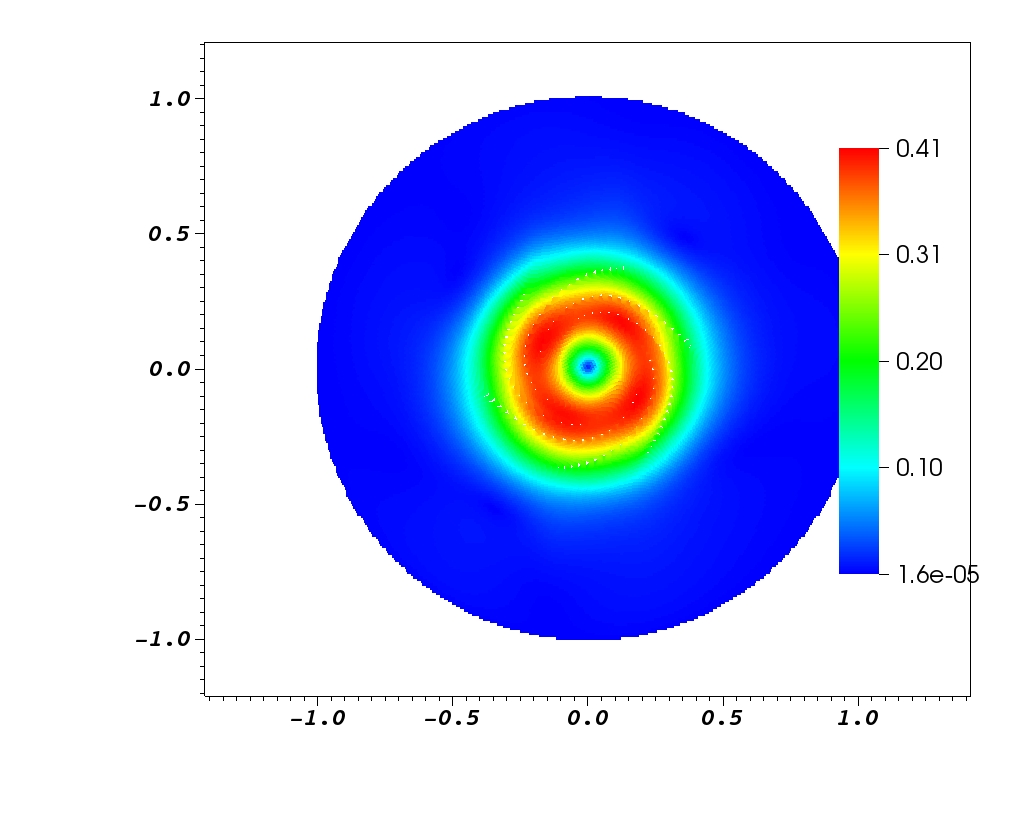}}
\caption{Exact solution (left column) and LP simulation  (right column) of the Gresho
  vortex. Top row images depict pressure and bottom row images depict
  density.}
\label{gresho_lp}
\end{figure}

\begin{figure}
\centering
\subfigure[MUSCL, pressure]{\label{gresho_lp_p0}\includegraphics[width=0.49\textwidth]{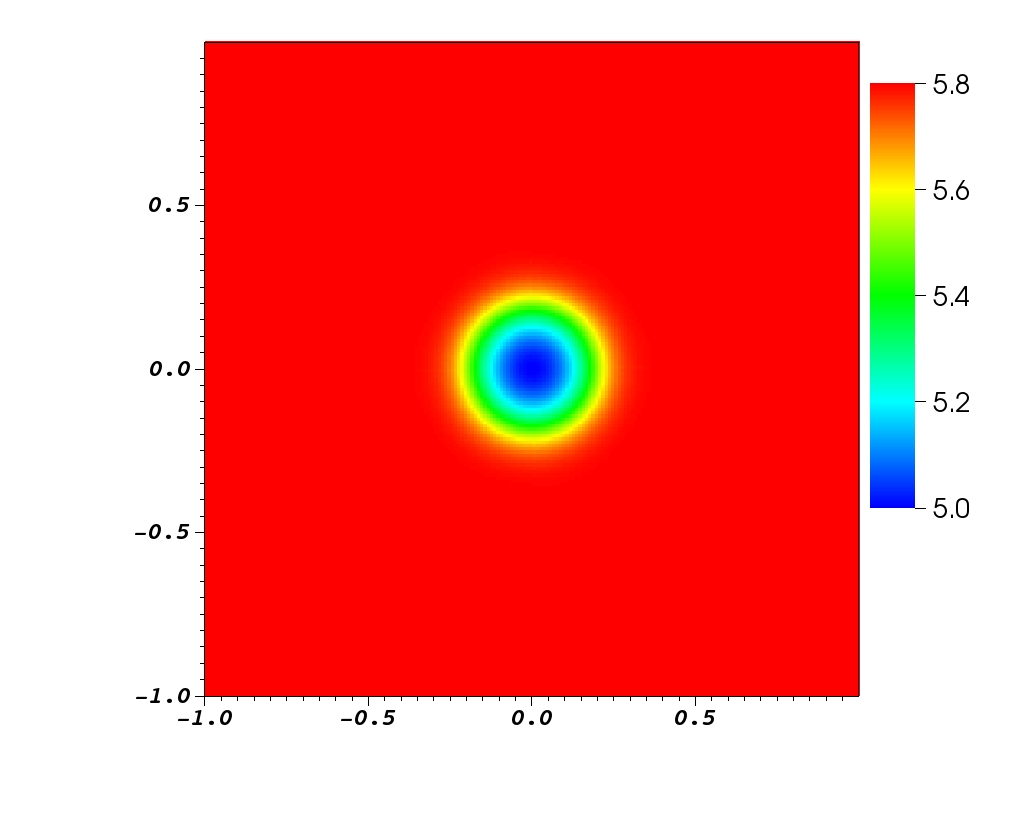}}
\subfigure[MUSCL, velocity]{\label{gresho_lp_p1}\includegraphics[width=0.49\textwidth]{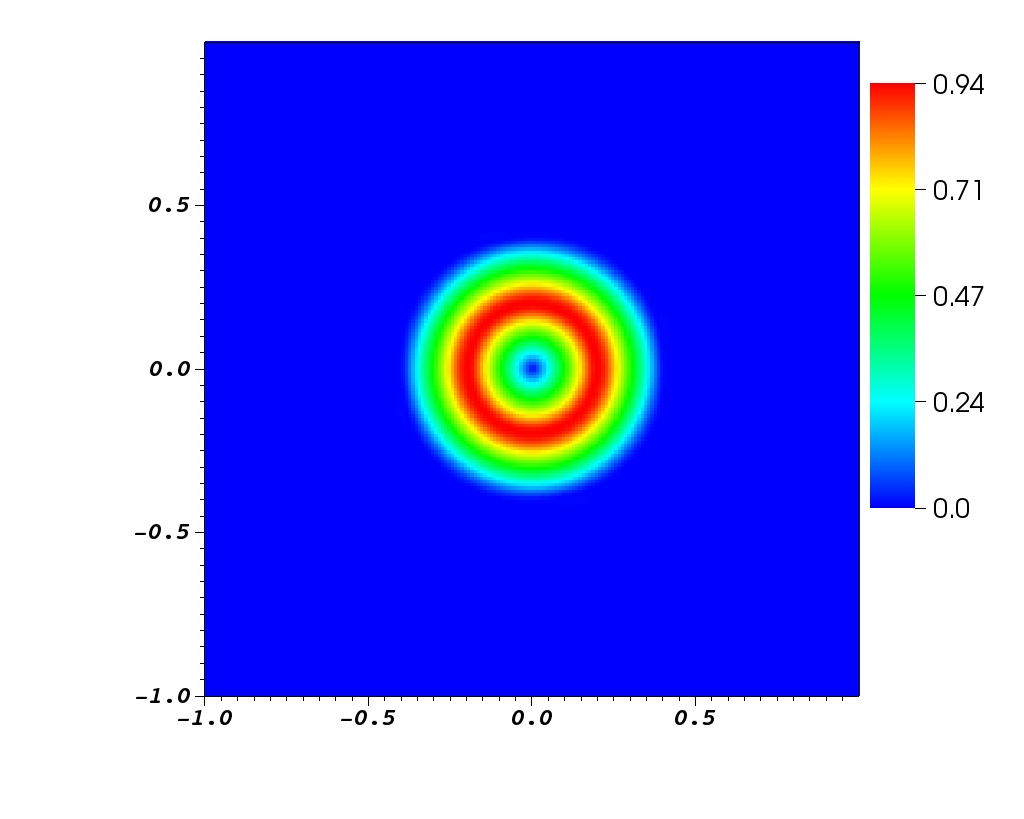}}
\subfigure[LF, pressure]{\label{gresho_lp_u0}\includegraphics[width=0.49\textwidth]{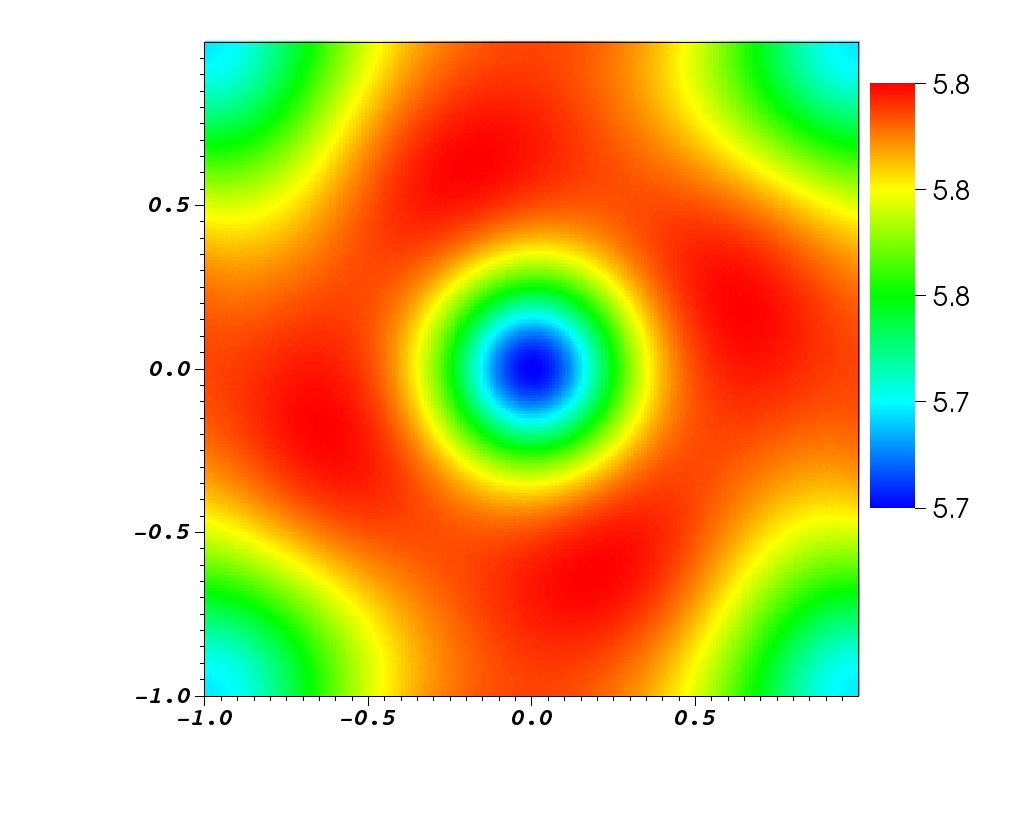}}
\subfigure[LF, velocity]{\label{gresho_lp_u1}\includegraphics[width=0.49\textwidth]{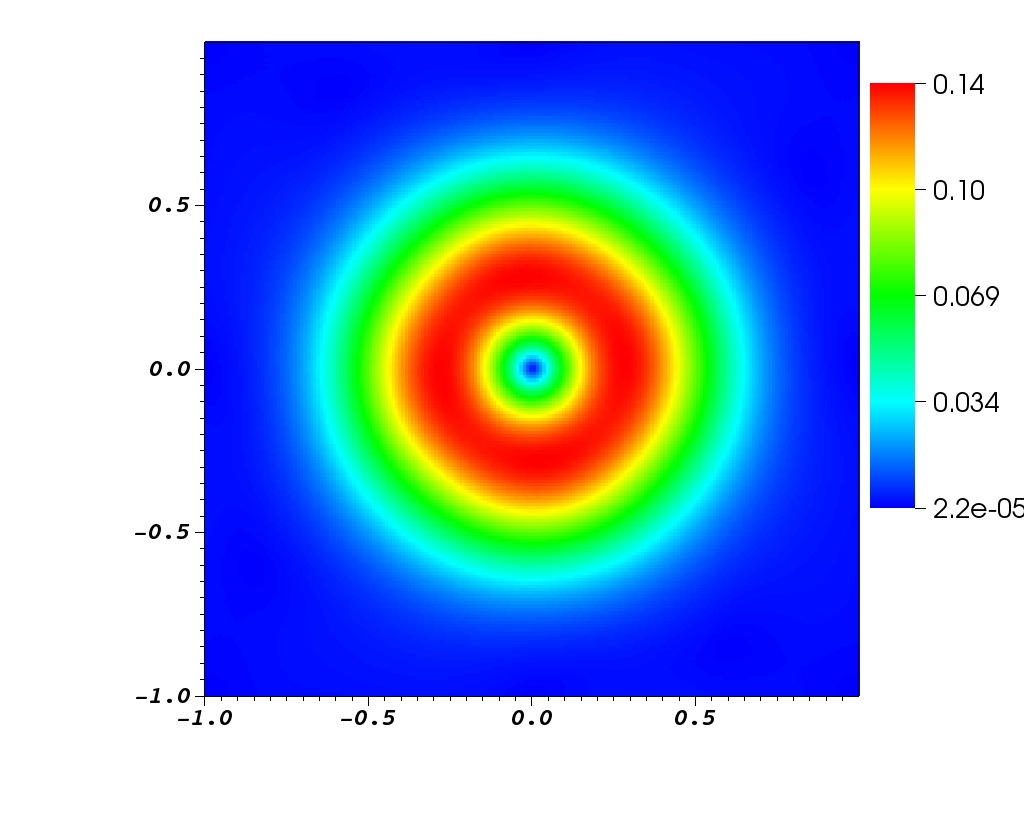}}
\caption{Simulations of the Gresho vortex using MUSC scheme (top row
  images) and Lax Friedrichs scheme (bottom row images).  Left column
  of images depict pressure and right column of images depict
  density.}
\label{gresho_grid}
\end{figure}

This problem, known to be notoriously difficult for SPH, tests the
accuracy of numerical scheme and its ability to preserve the symetry
and angular momentum. The empirical is order of accuracy of
second-order grid-based schemes for this problem, reported in
literature, is approximately 1.4.  Figure \ref{gresho_lp} shows the
analytic solution and numerical simulation results of the vortex after
one full rotation obtained with the Lagrangian particle method. The
convergence is of second order at the initial stages of rotation, and
it degrades to the first order in the later stages of rotation due to
particle redistribution.  For comparison, in Figure \ref{gresho_grid}
we show simulations performed using grid-based schemes at the same
numeircal resolution.

\subsection{2D Collision Between Two Circular Disks}
In previous test problems, particles are initialized using regular
distributions, such as the hexagonal packing, and slightly move with the flow.  
Nevertheless,  the magnitude of the particle
movement is quite restricted in previous tests, usually less than five percent of the
initial inter-particle-spacing.  In this section, a geometrically complex
two-dimensional problem with large particle movement and object
shape distortion is presented.  

The setup of the problem is as follows. Two fluid disks have initially uniform density $\rho = 1$ and
zero pressure, and  material properties described by the stiffened polytropic
 EOS with $\gamma=6$ and $P_{\infty}=7000$. The two disks move toward each other with the relative longitudinal velocity 
 of 20, but along lines that do not connect their centers. The time sequence shows the distortion 
 of disks after the collision. While no benchmark data exists for such a problem, we believe that the results are reasonable 
from physics point of view as they agree with theoretical estimates of achievable pressure peaks.
They demonstrate the ability of the proposed method to handle geometrically complex interfaces.
\begin{figure}
	\centering
	\includegraphics[width=0.32\textwidth]{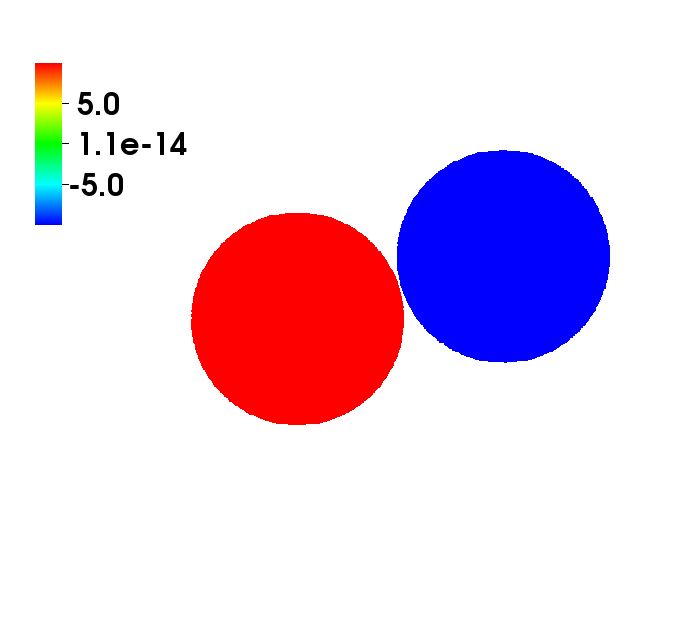} 
	\includegraphics[width=0.32\textwidth]{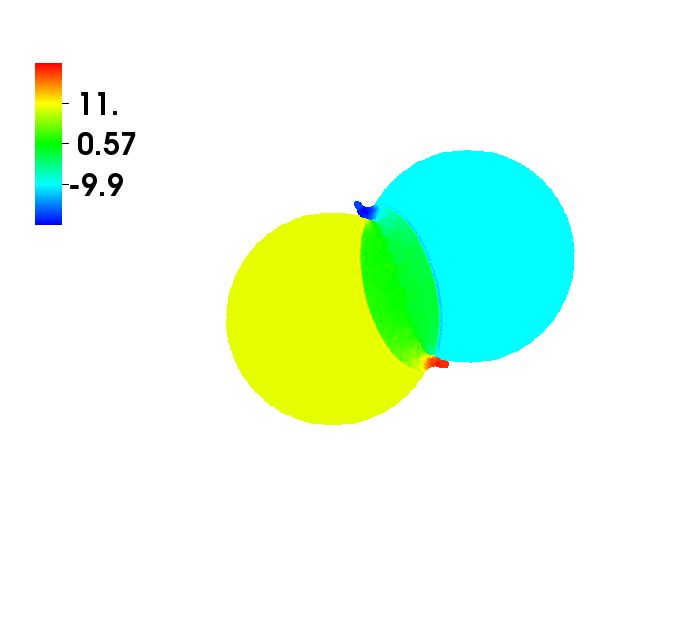} 
	\includegraphics[width=0.32\textwidth]{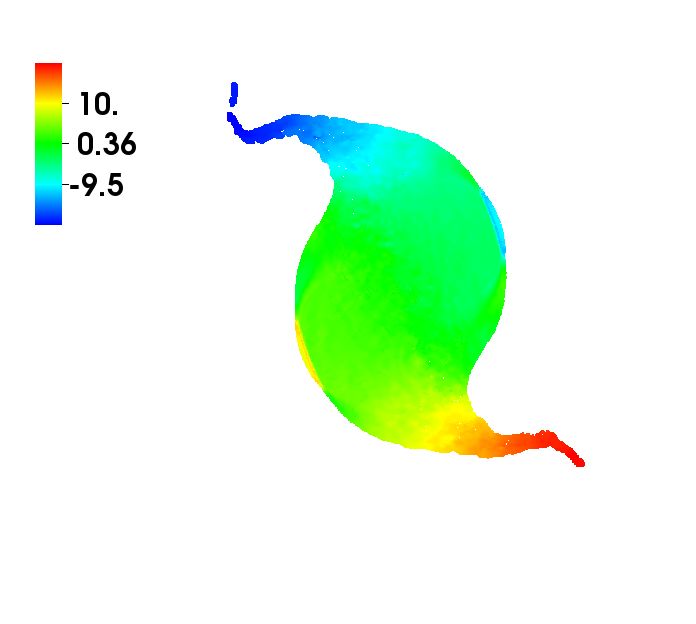} 
	\caption{2D simulation of collision of two disks. Velocity  distribution (10 m/s)
at initial time (left), 21  (middle), and 54  (right).}
	\label{fig:collide_spoly_free_2d_0}
\end{figure}

\subsection{Gas expansion into Vacuum}

\section{Conclusions and Future Work}

A Lagrangian particle method has been proposed for the simulation of Euler equations describing
compressible inviscid fluids or gases. By representing Lagrangian fluid
cells with particles, similarly to smoothed particle hydrodynamics, the method eliminates the mesh distortion problem of the
original Lagrangian method and is suitable for the simulation of
complex free surface flows. The main contributions of our method, which is different from SPH in all other aspects,
are (a) significant improvement of approximation of differential operators based on 
polynomial fits and the corresponding weighted least squares problems 
and convergence of prescribed order, (b) an upwinding second-order particle-based algorithm with limiter,
providing  accuracy and long term stability, (c)  elimination of the dependence on artificial parameters such as the
smoothening length in SPH, causing difficulties especially in the case
of large density changes, and (d) accurate resolution of states at free interfaces.
Numerical verification tests demonstrate the second convergence order of the method and its ability to resolve 
complex free surface flows. 

The Lagrangian particle method has numerous advantages compared to grid-based methods for the simulation of 
complex systems. It eliminates the need for complex and costly algorithms for the generation and adaptation of meshes,
provides continuos adaptivity to density changes, and is suitable for extremely non-uniform domains typical for 
astrophysics or high energy density applications. 
The algorithmic complexity of key particle methods
insignificantly increases with the increase of spatial dimensions, making a 3D code similar to a 1D code.
In addition, particle algorithms are independent of the geometric complexity of domains. In contrast, there
is a huge increase in algorithmic complexity of a 3D mesh generation and dynamic adaptation compared
to 1D as well as the increase associated with the geometric complexity of domains.

The future development
of the space-time discretization methods will explore new high resolution WENO-type solvers based on
irregularly placed particle nodes and symplectic integrators.
Our Lagrangian particle method is also generalizable to coupled multiphysics systems, including the
dynamics of plasmas, incompressible fluids, and fracture of solids. 

\vskip5mm
{\bf Acknowledgement.}
This research has been partially supported by the DOE Muon Accelerator Program.
This manuscript has been authored in part by Brookhaven Science Associates, LLC, under Contract No.
DE-SC0012704 with the US Department of Energy. The United States Government retains, and the
publisher, by accepting the article for publication, acknowledges, a world-wide license to publish or reproduce
the published form of this manuscript, or allow others to do so, for the United States Government
purpose.


\begin{thebibliography}{00}

\bibitem{RichtmyerMorton} R.D. Richtmyer, K.W. Morton, {\it Finite
  diference methods for initial value problems}, Interscience, New
  York - London - Sydney, 1967.

\bibitem{Lamb} H. Lamb, {\it Hydrodynamics}, Cambridge Univ. Press,

\bibitem{Hirt81} Hirt, C.W., Nichols, B.D. (1981), "Volume of fluid
  (VOF) method for the dynamics of free boundaries", J. Comput. Phys.,
  39 (1): 201–225

\bibitem{Osher02} S. Osher, R. Fedkiw, Level Set Methods and Dynamic
  Implicit Surfaces, 2002, Springer-Verlag.

\bibitem{ALE} C.W. Hirt, A.A. Amsden, J.L. Cook, An arbitrary
  Lagrangian-Eulerian computing method for a;; flow speeds,
  J. Comput. Phys, 135 (1997), 203-216.

\bibitem{FTcode} 
B.~Fix, J.~Glimm, X.~Li, Y.~Li, X.~Liu, R.~Samulyak, and Z.~Xu.
\newblock A {TSTT} integrated Frontier Code and its applications in
  computational fluid physics.
  \newblock {\em Journal of Physics: Conf. Series}, 16:471--475, 2005.

\bibitem{Monaghan92} J. J. Monaghan, Smoothed Particle Hydrodynamics,
  In: Annual review of astronomy and astrophysics. Vol. 30,
  p. 543-574. 1992.

\bibitem{Monaghan05} 
J. J. Monaghan, 
Smoothed particle hydrodynamics, Rep. Prog. Phys., 68 (2005), 1703 -- 1759.

\bibitem{Monaghan85} J. J. Monaghan, J. C. Lattanzio. A Refined
  Particle Method for Astrophysical Problems, Astronomy and
  Astrophysics, vol. 149, no. 1, Aug. 1985, p. 135-143.

\bibitem{Diltz99}
G. A. Dilts, Moving-least-squares particle hydrodynamics – I. 
Consistency and stability. 
Int. J. Num. Methods in Engineering, 44 (1999), 1115 -- 1155.

\bibitem{GIZMO}
P. F. Hopkins, GIZMO: a new class of accurate, mesh-free hydrodynamic 
simulation methods, Mon. Not, R. Astron. Soc., 2014.

\bibitem{GFD_ad}
F. Urena, J. J. Benito, L. Gavete, 
Application of the generalized finite difference method to solve the advection-diffusion equation,
J. Comp. Appl. Math., 235 (2011), 1849 - 1855.

\bibitem{Price_MHD} D. Price, Smoothed particle hydrodynamics and
  magnetohydrodynamics, arxiv.org:1012.1885v1

\bibitem{BeamWarming}
J.L. Steger, R.F. Warming,  
Flux vector splitting of the inviscid gas dynamic equations with application to finite-difference methods, 
J. Comput. Phys., 40 (1981), 263 - 293.

\bibitem{Strang_splitting}
G. Strang, On the construction and comparison of difference schemes,
SIAM J. Numerical Analysis,  5 (1968), 506-517.

\bibitem{BenitoUrena01} J.J. Benito, F. Urena, L. Gavete, Influence of
  several factors in the generalized finite difference method, Applied
  Math. Modeling, 25 (2001), 1039-1053.

\bibitem{CourantFriedrichs} R. Courant, K.O. Friedrichs, {\it
  Supersonic flow and shock waves}, Springer, 1999.

\bibitem{Fan93localpolynomial} Jianqing Fan, Theo Gasser, Irene
  Gijbels, Michael Brockmann, and Joachim Engel.  \newblock Local
  polynomial fitting: A standard for nonparametric regression, 1993.

\bibitem{Jiao08} X.~Jiao and H.~Zha.  \newblock Consistent computation
  of first- and second-order differential quantities for surface
  meshes.  \newblock {\em ACM Solid and Physical Modeling Symposium},
  2008. (Available at arxiv.org/abs/0803.2331)

\bibitem{GolubVanLoan96} G.~GOLUB and C.~F.~VAN LOAN.  \newblock {\em
  Matrix Computations}.  \newblock Johns Hopkins University Press,
  Baltimore, MD, 1996.

\bibitem{Samet_octree}
H. Samet, The design and analysis of spatial data structures,
Adison-Wesley Publishing Company, 1990.

\bibitem{Shu_WENO} Chi-Wang Shu, Essentially Non-Oscillatory and
  Weighted Essentially Non-Oscillatory Schemes for Hyperbolic
  Conservation Laws, NASA/CR-97-206253, ICASE Report No. 97-65,
  Nov. 1997

\bibitem{LiskaWen03} 
R. Liska, B. Wendroff, 
Comparison of several difference schemes on 1D and 2D test problems for the Euler equations,
SIAM J. Sci. Comput., 25 (2003), 995–1017. 



\end{thebibliography}
\end{document}